\documentclass[11pt]{article}
\usepackage{amscd,amssymb,amsmath,verbatim,amsthm}
\usepackage{times}
\usepackage{enumerate}
\usepackage[11pt,curve,matrix,arrow,frame,graph]{xy}
\newcounter{intro}

\newtheorem{thm}{Theorem}[section]
\newtheorem{lem}[thm]{Lemma}
\newtheorem{prop}[thm]{Proposition}

\newtheorem{rem}[thm]{Remark}

\numberwithin{equation}{section}   

\newcounter{counteroman}

\newcommand{\cref}[1]{Corollary~\ref{#1}}

\newcommand{\rref}[1]{Remark~\ref{#1}}

\newcommand{\R}{\mathbb{R}}

\newcommand{\bB}{\mathbb{B}}

\newcommand{\bpm}{\begin{pmatrix}}
\newcommand{\epm}{\end{pmatrix}}

\let\eps=\varepsilon

\DeclareMathOperator{\scal}{Scal}
\DeclareMathOperator{\diam}{diam}
\DeclareMathOperator{\Lip}{Lip}

\numberwithin{equation}{section}

\renewcommand{\tilde}{\widetilde}

\renewcommand{\hat}[1]{\widehat{#1}}

\newcommand\paperintro%
        {%
         }
\newcommand\paperbody%
        {%
         }

\newcommand\cC{\mathcal{C}}
\newcommand\cD{\mathcal{D}}
\newcommand\cE{\mathcal{E}}

\newcommand\cH{\mathcal{H}}

\newcommand\cL{\mathcal{L}}

\newcommand\cU{\mathcal{U}}

\DeclareMathAlphabet{\mathpzc}{OT1}{pzc}{m}{it}

\newcommand{\RR}{\mathbb R}

\newcommand{\calC}{{\mathcal C}}

\newcommand{\calE}{{\mathcal E}}

\begin{document}
\title{The Yamabe problem on Dirichlet spaces}
\author{Kazuo Akutagawa \thanks{akutagawa@math.titech.ac.jp}\\ Tokyo Institute of Technology \and 
Gilles Carron \thanks{Gilles.Carron@math.univ-nantes.fr} \\ Universit\'e de Nantes  \and 
Rafe Mazzeo \thanks{mazzeo@math.stanford.edu}\\ Stanford University}

\date{May, 2013}




\maketitle


\section{Introduction} 
This article is a successor to our previous paper \cite{ACM} and continues the theme of generalizing the 
Yamabe problem to various classes of singular spaces. In that earlier paper we considered this problem on 
`almost smooth' metric-measure spaces which satisfy a small set of additional structural hypotheses. 
As part of this, we defined the local Yamabe invariant $Y_{\ell}(M,[g])$, which is a generalization of the
quantity $Y(S^n)$ which plays a key role in the standard Yamabe problem, and then established solvability 
of the Yamabe problem for any metric $g$ on the smooth locus of one of these spaces provided it satisfies 
$-\infty < Y(M, [g]) < Y_{\ell}(M, [g])$. As the main application there, we find Yamabe minimizers on certain
stratified spaces with iterated edge metrics.

In the present article we consider this problem in a more general setting, on the class of {\it Dirichlet spaces} which 
satisfy a few additional structural properties. Our main results here again concern the {\it generalized Aubin inequality}, 
in particular its role in establishing {\it existence} of minimizing solutions for the Yamabe energy, and we also consider 
the {\it regularity} for (not necessarily minimizing) critical points of this energy. 

Let us begin by recalling the standard Yamabe problem. Consider the functional 
\[
E(g) := \frac{\int_M {\rm Scal}_g \, d\mu_g}{\textrm{Vol}_g(M)^{(n-2)/n}} 
\]
on the space $\mathcal{M}(M)$ of all Riemannian metrics on the compact (smooth) manifold $M^n$ with $n \geq 3$. 
Here, ${\rm Scal}_g, d\mu_g$ and 
$\textrm{Vol}_g(M)$ are the scalar curvature, volume form and volume of $(M, g)$. This is called
the (normalized) Einstein-Hilbert functional, and its critical points are the Einstein metrics on $M$. 

This functional is unbounded both above and below, so it is reasonable to search for critical points using 
a max-min scheme. Consider the quantity
\[
Y(M, C) := \inf_{\tilde{g} \in C} E(\tilde{g}),
\]
the infimum of $E$ on any conformal class $C = [g] := \{ \textrm{e}^{2f}\cdot g\ |\ f \in \calC^{\infty}(M) \}$.
This is called the {\it Yamabe invariant} (or {\it Yamabe constant} or {\it conformal Yamabe invariant}) of $C$.
We then define
\[
Y(M) := \sup_{C \in \mathcal{C}(M)} \inf_{g \in C} E(g) = \sup_{C \in \calC(M)} Y(M, C)\,,  
\]
where $\calC(M)$ is the space of all conformal classes on $M$. This is called the {\it Yamabe invariant} (or 
{\it $\sigma$-invariant} or {\it smooth Yamabe invariant}) of $M$, see \cite{Ko}, \cite{Sc-2}.  

The {\it Yamabe problem} concerns the first part of this, namely whether it is possible to find
a metric which minimizes $E$ in a given conformal class $C$. Such a metric has constant
scalar curvature, and conversely, any constant scalar curvature metric is at least a critical
point for $E$ in its conformal class. The second step, showing that one can find a metric $g$
which attains the max-min, so that $g$ is Einstein (and $E(g) = Y(M)$) is significantly more
difficult. We refer to \cite{LeBrun}, \cite{Anderson}, \cite{ADH} for some significant progress here. 

It is now well known, through successive work of Yamabe, Trudinger, Aubin and Schoen,
see \cite{LP}, \cite{Aubin-Book} for details, 
that each conformal class $C$ contains a minimizer $\hat{g}$ of $E$ restricted to that conformal class, 
called the {\it Yamabe metric} of that class, and 
\[
{\rm Scal}_{\hat{g}} = Y(M, C)\cdot \textrm{Vol}_{\hat{g}}(M)^{-2/n}\,. 
\]

When studying sequences of Yamabe metrics $g_j$ satisfying certain geometric non-collapsing assumptions, and 
with $E(g_j) \to Y(M)$, one is led to consider limit spaces which are {\it Riemannian orbifolds} (or {\it Riemannian multi-folds}, 
{\it manifolds with conic singularities}, {\it simple edge spaces}, and more general {\it iterated edge spaces}). 
This makes it natural to study the Yamabe problem directly on these and more general singular spaces, cf.\ 
\cite{Ak-2}, \cite{AB}, \cite{ACM}, \cite{JR}, \cite{Vi}.

In our previous paper~\cite{ACM} we consider the Yamabe problem on a compact metric-measure space $(M, d,\mu)$ 
which has a compatible smooth Riemannian metric $g$ on an open dense subset; we call this an {\it almost smooth} 
metric-measure space. Assuming also that this space is Ahlfors $n$-regular, satisfies a Sobolev inequality,
and with certain growth conditions on $\mathrm{Scal}_g$, but without specific information about the singular set 
of $M$, we define the {\it local Yamabe invariant} $Y_{\ell}(M, d, \mu)$. Roughly speaking,
this is the infimum of the Yamabe invariants of each of the tangent cones to $M$. 
When $M$ is smooth, $Y_\ell(M,[g])$ equals the Euclidean Yamabe invariant $Y(\RR^n)$, or equivalently,
the Yamabe invariant $Y(S^n, [g_0])$ of the round sphere $(S^n, g_0)$. {\it Aubin's inequality} \cite{Aubin} states that
$Y(S^n, [g_0])$ is the supremum of the set of values of the Yamabe invariants over all compact smooth conformal $n$-manifolds: 
\[
Y(M, C) \leq Y(S^n, [g_0])\quad \mbox{for every}\ \ (M, C)\,. 
\]
As in the smooth case, we can define the Yamabe invariant $Y(M, d, \mu)$ of a compact metric-measure space $(M, d, \mu)$,
and it is not hard to show that the analogous Aubin-type inequality
\begin{equation}
Y(M, d, \mu) \leq Y_{\ell}(M, d, \mu)
\label{Aubinineq}
\end{equation}
is still valid. This local Yamabe invariant contains much information about the metric near the singular points of $M$.
In \cite{ACM}, we showed that if $(M, d, \mu)$ is almost smooth and satisfies the extra conditions noted above, 
and if $Y(M, d, \mu) < Y_{\ell}(M, d, \mu)$, then the energy $E$ attains its minimum in that conformal class.
We also proved that solutions of the Yamabe equation on $(M, d, \mu)$ are bounded and uniformly positive.

We generalize this yet further here and consider a Yamabe-type problem on a so-called Dirichlet space $(M, \mu, \cE)$, i.e.\ 
a finite measure space $(M, \mu)$ equipped with a Dirichlet form $\cE$ on $L^2(M, \mu)$, with the scalar curvature replaced 
by a potential $V$. Assuming a few other conditions on the space and potential, we define a Yamabe invariant $Y(V)$ of 
$(M, \mu, \cE, V)$, and then consider the corresponding Yamabe-type problem. After proving the generalization of 
\eqref{Aubinineq}, we show that if this inequality is strict, then (again under certain additional assumptions), this 
generalized Yamabe problem admits a minimizer. We also prove the boundedness, uniform positivity and H\"older continuity 
of more general solutions of the associated Yamabe equation. 

This paper is organized as follows: \S 2 reviews the necessary terminology and defines the generalized Yamabe problem 
on a Dirichlet space; in \S 3 we establish the Aubin inequality and prove existence of minimizers of the generalized Yamabe problem;
\S 4 contains proofs of the regularity results for solutions of the Yamabe type equation; finally, in \S 4, we 
present some examples of this generalized Yamabe problem. 

\noindent 
{\bf Acknowledgements.} 
The authors have been supported by the following grants: K.A. through the Grant in-Aid for Scientific Research 
(B), JSPS, No. 24340008; G.C. through the ANR grant ACG: ANR-10-BLAN 0105; R.M. by the NSF under
DMS-1105050.

\section{A generalized Yamabe problem} 
We begin by presenting some terminology which allows us to pose the generalized Yamabe 
problem on a Dirichlet space.

\subsection{Dirichlet spaces} 
We first review some classical facts about Dirichlet spaces; \cite{fuku} is a comprehensive reference for 
this material, but see also \cite[page 209]{RS4} or \cite{Besson}, which is sufficient for what we do here.

Let $(M, \mu)$ be a finite measure space, and consider a nonnegative closed symmetric bilinear form $\cE$ 
defined on a dense subspace $\cD(\cE) \subset L^2(M, \mu)$; thus $\cE : \cD(\cE) \times \cD(\cE) \rightarrow \overline{\RR^+}$.
We refer to this simply as a {\it closed symmetric form} on $L^2(M, \mu)$,  and identify $\cE$ with the 
corresponding quadratic form $\cE(\varphi, \varphi)$.  Because this quadratic form is semibounded, the Friedrichs 
extension procedure determines a selfadjoint operator $L\colon\cD(L)\rightarrow L^2(M,\mu)$, 
with domain $\cD(L)$ consisting of all functions $v\in \cD(\cE)$ such that 
\[
|\cE(v,\varphi)|\le C \|\varphi\|_{L^2}\quad \mbox{for every}\ \ \varphi\in \cD(\cE).
\]
for some constant $C$ (depending on $v$ and $\calE$, but not $\varphi$). This is the generator of $\cE$. 
A closed symmetric form on $L^2(M, \mu)$ is called a {\it Dirichlet form} if its generator $L$ is {\it subMarkovian}, i.e.\ provided
the semigroup $e^{-tL}$ satisfies 
\[
0\le v\le 1\Rightarrow 0\le e^{-tL}v\le 1\,\,.
\]
According to the Beurling-Deny criteria, this is equivalent to the following: 
\begin{enumerate}[i)]
\item $v\in \cD(\cE)\Rightarrow |v|\in \cD(\cE)\ \mathrm{and}\ \cE(|v|)\le \cE(v)$
\item $v\in\cD(\cE)\ \mathrm{and}\ v \ge 0\Rightarrow v_1 := \inf\{v, 1\} \in \cD(\cE) \ \mathrm{and}\ \cE(v_1)\le \cE(v)$. 
\end{enumerate}

A triple $(M,\mu,\cE)$ with all these properties is called a {\it Dirichlet space}.  

\subsection{The Sobolev inequality}
Suppose that $(M,\mu,\cE)$ is a Dirichlet space for which a Sobolev inequality holds. This means that there exist $\nu > 2$ 
and $A, B > 0$ such that 
\begin{equation}\label{eq:Sobolev}
A \|v\|^2_{L^{\frac{2\nu}{\nu-2}}}\le A\cE(v)+B\|v\|_{L^2}^2\quad {\rm for}\ \  v \in \cD(\cE).
\end{equation}
Following Nash \cite{Nash},  the heat semi-group $\left\{e^{-tL}\right\}_{t\ge 0}$ then necessarily satisfies an ultracontractive estimate:
there exists a constant $C$ such 
\begin{equation}\label{contractive}
\left\|e^{-tL}\right\|_{L^1\to L^\infty}\le \frac{C}{t^{\nu/2}}\, ,\quad 0 < t < 1. 
\end{equation}
It is now known, through work of Varopoulos \cite{Varopoulos}, that \eqref{eq:Sobolev} and \eqref{contractive} are in fact equivalent.

It is straightforward to show that \eqref{eq:Sobolev} implies the following compactness result, see \cite[Proposition 1.6]{ACM}: 
\begin{prop} If \eqref{eq:Sobolev} holds, then the inclusion
\[
\cD(\cE) \longrightarrow L^{\frac{2p}{p-2}}(M)
\]
is compact for any $p \in (\nu ,\infty]$.
\end{prop}

\subsection{Schr\"odinger operators} 
A nonnegative measurable function $W$ is said to be \textit{relatively form bounded} with respect to $\cE$ if there exists 
some constant $D>0$ such  that 
\[
\int_M Wv^2d\mu\le D( \cE(v)+ ||v||_{L^2}^2) \quad {\rm for}\ \  v \in \cD(\cE);
\]
similarly, $W$ is \textit{infinitesimally form bounded} with respect to $\cE$ if for any $\eps>0$ there exists $c(\eps)$ such that :
\[
\int_M Wv^2d\mu\le \eps\,  \cE(v)+ c(\eps) \int_M v^2d\mu\quad {\rm for}\ \ v \in \cD(\cE).
\]

Since $\cD(\cE) \hookrightarrow L^2(M)$ is compact, $W$ is infinitesimally form bounded with respect to $\cE$ if 
and only if the operator $(L+1)^{-\frac12}W(L+1)^{-\frac12}$ is compact on $L^2$.

If $V$ is a real-valued integrable function on $M$ and its nonpositive part $V_- := \sup\{0, -V\}$ is relatively form 
bounded with respect to $\cE$, we define the quadratic form
\[
\cE_V(v)=\cE(v)+\int_M Vv^2d\mu
\]
on the domain $\cD(\cE_V)=\{ v \in \cD(\cE) : \int_M Vv^2\, d\mu<\infty \}$.  As before, $\cE_V$ is densely defined, closed
and semibounded, so we can define the self-adjoint operator $L+V$  by the Friedrichs procedure. 
 
\subsection{The generalized Yamabe problem}
Let $V$ be integrable and suppose that $V_-$ is relatively form bounded; suppose too that the Dirichlet space $(M, \mu, \cE)$ 
satisfies the Sobolev inequality~\eqref{eq:Sobolev}. We then define the {\it Yamabe invariant} associate to the operator $L+V$ :
\[
Y(V)=\inf\left\{ \cE_V(v) : v \in \cD(\cE)\,\, \mathrm{and} \,\, \|v\|_{L^{\frac{2\nu}{\nu-2}}} = 1 \right\}\,.
\]
Note that \eqref{eq:Sobolev} implies immediately that
\[
Y(V)\ge -D\max\{A,B\}\,\,.
\]

We wish to whether there exists $u\in \cD(\cE_V)$ such that
\[
\cE_V(u)=Y(V)\quad {\rm and}\quad \|u\|_{L^{\frac{2\nu}{\nu-2}}}=1\,\,.
\]
Since
\[
\cE_V(|u|)\le \cE_V(u),
\]
we can always assume that any such minimizer must be nonnegative. This minimizer must satisfy
the Euler-Lagrange equation
\begin{equation}
\label{Yamabe}
\cE_V(u,\varphi)=Y(V)\int_M u^{\frac{\nu+2}{\nu-2}} \varphi\, d\mu\quad \mbox{for every}\ \ \varphi\in \cD(\cE)\,\,.
\end{equation}
Note that by the Sobolev and H\"older inequalities, the right hand side is finite.

\section{Existence of minimizers} 
\subsection{Existence theorem} 
\begin{thm} Let $(M,\mu,\cE)$ be a Dirichlet space with Sobolev inequality~$\eqref{eq:Sobolev}$ for some $\nu>2$ 
and positive constants $A, B$. Let $V$ be an integrable function whose nonpositive part $V_-$ is infinitesimally form 
bounded with respect to $\cE$. Assume that 
\begin{equation}
Y(V) < \frac{1}{A}\,\,. 
\label{strAub}
\end{equation}
Then there exists $v\in \cD(\cE)$ such that 
\[
\cE_V(v)=Y(V)\quad \mathrm{and}\quad \|v\|_{L^{\frac{2\nu}{\nu-2}}}=1\,\,.
\]
\end{thm}
\begin{rem}\label{Optconstant} 
The hypothesis \eqref{strAub} can be rephrased in terms of the {\rm optimal Sobolev constant} $A_{\mathrm{opt}}$. This is,
by definition, the smallest constant such that for every $A > A_{\mathrm{opt}}$ there exists $B > 0$ such that the
Sobolev inequality $\eqref{eq:Sobolev}$ holds with that choice of $A$ and $B$. 
We also write
\[
A_{\mathrm{opt}}=\frac{1}{\alpha(\cE)}\,\,,
\]
where
\[
\alpha(\cE)=\lim_{t\to \infty} Y(t)\,,
\]
i.e.\ the limit of the Yamabe invariants associated to the constant potentials $V \equiv t$.
Another characterization is that
\[
A_{\mathrm{opt}}=\lim_{t\to +\infty}\left\| \left(\sqrt{L}+t\right)^{-\frac12}\right\|^2_{L^2\to L^{\frac{2\nu}{\nu-2}}}.
\]
\end{rem}
\proof
Using the infinitesimal form boundedness of $V_-$,  we see that if $\hat A > A$ then there exists a positive 
constant $\hat B$ such that :
\begin{equation}\label{eq:SobV}
\|v\|^2_{L^{\frac{2\nu}{\nu-2}}}\le \hat A \cE_V(v)+\hat B\|v\|_{L^2}^2\quad \mbox{for all}\ \  v \in \cD(\cE)\,\,.
\end{equation}
Choose $\hat A>A$ so that $\hat A Y(V)<1$. Since the embedding $\cD(\cE_V)\rightarrow L^2$ is compact, 
we can find a minimizing sequence $u_\ell \in \cD(\cE_V)$ and $u\in \cD(\cE_V)$ such that 
\begin{itemize}
\item[a)] $u_\ell \rightharpoonup u$ weakly in $\cD(\cE_V)$;
\item[b)] $u_\ell \to u$ strongly in $L^2$;
\item[c)] $u_\ell \to u$ a.e.;
\item[d)] $\|u_\ell\|_{L^ {\frac{2\nu}{\nu-2}}}=1$.
\end{itemize}
By d), 
\[
\cE_V(u_\ell-u)=\cE_V(u_\ell)-\cE_V(u)+\eps_\ell, \qquad \mbox{where}\quad \lim_{\ell \to \infty} \eps_\ell=0.
\]

We now appeal to a very useful result of Brezis and Lieb \cite{BL} (we are grateful to E. Hebey for 
pointing us to this), which gives
\[
\lim_\ell \left( \|u_\ell\|^ {\frac{2\nu}{\nu-2}}_{L^ {\frac{2\nu}{\nu-2}}} - 
\|u_\ell-u\|^ {\frac{2\nu}{\nu-2}}_{L^ {\frac{2\nu}{\nu-2}}} \right) =\|u\|^ {\frac{2\nu}{\nu-2}}_{L^ {\frac{2\nu}{\nu-2}}}\,.
\]
Hence, setting $I=\|u\|^ {\frac{2\nu}{\nu-2}}_{L^ {\frac{2\nu}{\nu-2}}}$, then 
\[
\lim_\ell \|u_\ell-u\|^ {\frac{2\nu}{\nu-2}}_{L^ {\frac{2\nu}{\nu-2}}}=1-I\,.
\]
Now apply the Sobolev inequality \eqref{eq:SobV} to $u_\ell-u$ and pass to the limit $\ell \to \infty$ to get
\[
(1-I)^{1-\frac2\nu}\le \hat A Y(V) -\hat A \cE_V(u)\,.
\]
On the other hand, by definition,
\[\cE_V(u)\ge Y(V) I^{1-\frac2\nu}\,,
\]
so putting these together and recalling the choice of $\hat A$ yields
\[
(1-I)^{1-\frac2\nu}+\hat A Y(V) I^{1-\frac2\nu}\le  \hat A Y(V) < 1.
\]
This forces $I=1$, hence $u \not \equiv 0$, and since $Y(V)\ge\cE_V(u)$, we conclude that $u$ is a minimizer 
for $\calE_V$. 
\endproof
\subsection{On the optimal Sobolev constant}
We now turn to a more careful discussion of the optimal Sobolev constant $A_{\mathrm{opt}}$ introduced in \rref{Optconstant}. 
We assume henceforth that $M$ is a \emph{compact topological space} and $\mu$ is a {\it Radon measure}, and moreover, 
that the Dirichlet space is \emph{regular} and \emph{strongly local}. These last two conditions are:
\begin{itemize}
\item (Regularity) $\cD(\cE)\, \cap \, \cC^0(M)$ is dense in both $\cD(\cE)$ with $\cE_1$-norm and $\cC^0(M)$ with uniform norm;
\item (Strong locality) if $u,v\in \cD(\cE)$  and if $u$ is constant in a neighborhood of $\mathrm{supp}(v)$, then $\cE(u,v)=0$. 
\end{itemize}
These conditions guarantee the existence of a bilinear form $d\gamma$, the so-called 
{\it the energy measure}, from $\cD(\cE)\times \cD(\cE)$ to the set of Radon measures on $M$, such that 
\[
\cE(u,v)=\int_M d\gamma(u,v)\quad {\rm for}\ \  u, v \in \cD(\cE)\,.
\]
If the energy measure is absolutely continuous with respect to $d\mu$, Bakry and Emery \cite{BE} call 
this bilinear form the {\it carr\'e du champ}. The energy measure is determined by the identity
\[
\cE(\phi u,u)-\frac12 \cE(\phi,u^2)=\int_M \phi \,d\gamma(u,u)\quad {\rm for}\ \ u\in \cD(\cE)\ \mbox{and}\ 
\phi\in \cD(\cE) \cap \ \cC^0(M)\,.
\]
The energy measure  satisfies the Leibniz and chain rules: 
\[
\begin{split}
d\gamma(uv,w) & =ud\gamma(v,w)+vd\gamma(u,w)\quad {\rm for}\ \  u, v, w\in \cD(\cE)\, \\[3mm]
d\gamma(f(u),v) & =f'(u)d\gamma(u,v)\quad {\rm for}\ \ u, v \in \cD(\cE),\ \mbox{and}\ f\in \Lip(\R)\,.
\end{split}
\]

A regular, strongly local Dirichlet space $(M,\mu,\cE)$ has an intrinsic pseudo-distance defined by
\[
d(x,y) = \sup\left\{u(x)-u(y) : u\in \cD(\cE)\cap \ \cC^0(M)\ \ \mathrm{and}\ \ d\gamma(u,u)\le d\mu\right\};
\]
the comparison $d\gamma(u,u)\le d\mu$ here means that there exists a function $f \le 1$ such that $d\gamma(u,u)=fd\mu$. 

If this pseudo-distance is compatible with the topology of $M$, then for any $y\in M$, the function $r_y = d(y, \cdot )$ 
satisfies $d\gamma(r_y,r_y) \le d\mu$ \cite{SturmI}. If $U$ is open in $M$, we define 
\[
\begin{split} 
S(U) & =\inf\left\{\cE(u): \|u\|_{L^ {\frac{2\nu}{\nu-2}}}=1\ \ {\rm and}\ \ \mathrm{supp}\, u\subset U\right\}\,, \\
Y(U) & =\inf\left\{\cE_V(u): \|u\|_{L^ {\frac{2\nu}{\nu-2}}}=1\ \ {\rm and}\ \ \mathrm{supp}\, u\subset U\right\}\,.
\end{split}
\]
We now adapt the proof of \cite[Proposition 1.4]{ACM}, using cutoffs of these distance functions, to obtain
\begin{prop} Let $(M,\mu,\cE)$ be a regular, strongly local Dirichlet form with intrinsic distance compatible 
with the topology of $M$. Then 
\[
A_{\mathrm{opt}}=\sup_{x\in M} \lim_{r \searrow 0}A_{\mathrm{opt}}(B(x,r))\,,
\]
where $B(x,r)$ denotes the metric ball of radius $r$ centered at $p$. 
If $A_{\mathrm{opt}}$ is finite, then 
\[
A_{\mathrm{opt}}=\frac{1}{S_\ell}\,,\quad {\rm where}\ \ S_{\ell} := \inf_{x\in M} \lim_{r \searrow 0} S(B(x,r))\,.
\]
Moreover, if $|V|$ is infinitesimally form bounded with respect to $\cE$, then 
\[
S_{\ell} = Y_{\ell} := \inf_{x\in M} \lim_{r \searrow 0} Y(B(x,r))\,.
\]
\end{prop}

\section{Regularity of solutions}  
In this section, we now prove various facts about regularity of solutions (which are not necessarily minimizers) of this generalized 
Yamabe equation. Note that this equation can be rewritten as
\begin{equation}\label{Yamaberw}
Lu = Wu,\quad {\rm where}\ \ W=-V+Y(V) u^{\frac{4}{\nu-2}}\,. 
\end{equation} 
Some of our results will follow from regularity results for solutions of this linear equation. 
\subsection{Boundedness} 
\subsubsection{General results} 
\begin{prop} 
Let $(M,\mu,\cE)$ be a Dirichlet space with Sobolev inequality \eqref{eq:Sobolev}. 
Let $W$ be a nonnegative measurable function with $W\in L^q$ for some $q>\nu/2$. Assume that $u\in \cD(\cE)$ is a nonnegative 
function satisfying 
\begin{equation}\label{subL} 
Lu \le Wu\,. 
\end{equation}
Then $u \in L^\infty$, and moreover, 
\[
\|u\|_{\infty} \leq C \|u\|_{2},
\]
where the constant $C$ depends only on $\|W\|_{L^q}$, $n, q$ and the constants $A, B$. 
\end{prop}
This follows from the Gagliardo-Nirenberg inequality \cite{Coulhon}; the proof is in \cite{carher}, but
for the sake of completeness, we sketch the proof here as well.
\proof 
The inequality \eqref{subL} means that for all nonnegative $\varphi\in \cD(\cE)$, 
\[
q(u,\varphi)\le \int_M W u\varphi \, d\mu\,.
\]
The Sobolev inequality implies that 
\begin{equation*}
\left\|e^{-tL}\right\|_{L^1\to L^\infty}\le \frac{C}{t^{\nu/2}}\quad {\rm for}\ \ t \in(0,1)\,,
\end{equation*}
and hence by interpolation, if $1 \le r < s $, then 
\begin{equation}\label{contractive2}
\left\| e^{-tL}\right\|_{L^r\to L^s}\le \frac{C}{t^{\frac{\nu}{2}\left(\frac{1}{r}-\frac{1}{s}\right)}}\quad {\rm for}\ \ t \in (0,1)\,. 
\end{equation}

Clearly, 
\[
Le^{-tL}u=e^{-tL}Lu\le e^{-tL}Wu\,,
\]
and hence  
\begin{equation*}
u=e^{-L}u+\int_0^1 e^{-tL}Lu\, dt \le e^{-L}u+\int_0^1 e^{-tL}Wu\, dt .
\end{equation*}
Now introduce
\[
T(f)=e^{-L}f+\int_0^1 e^{-tL}Wf\, dt .
\]
From \eqref{contractive} and \eqref{contractive2}, it follows that if 
$f\in L^s$ with $\frac1s<\frac{2}{\nu}-\frac1q$, then $Tf\in L^\infty$. 
Indeed, if $r$ is determined by $r^{-1}=s^{-1}+q^{-1}$, then 
\begin{equation*}\begin{split}
\|Tf\|_{L^\infty}&\le C \|f\|_{L^s}+\int_0^1\left\|e^{-tL}\right\|_{L^r\to L^\infty}\,\|W\|_{L^q}\,\|f\|_{L^s}\,dt\\
&\le  C \|f\|_{L^s}+\int_0^1C t^{-\frac{\nu}{2}\left(\frac1q+\frac1s\right)}\,\|W\|_{L^q}\,\|f\|_{L^s}\, dt\\
&\le C (1+\|W\|_{L^q})\|f\|_{L^s}\,.
\end{split}\end{equation*}

A similar argument shows that if 
\[
f\in L^s\quad \mathrm{with}\ \ \frac1s>\frac{2}{\nu}-\frac1q\,,
\]
then
\[
Tf \in L^r\quad {\rm for}\ \ r \ge 1\,\, \mathrm{and}\,\, \frac1r > \frac1s + \frac1q - \frac{2}{\nu}\,.
\]
Hence, from $u \in L^2$, we obtain that $u \in L^\infty$ in a finite number of steps.
\endproof

\begin{rem}\label{nu2}
It is easy to show using \eqref{eq:Sobolev} that if $W\in L^{\frac{\nu}{2}}$, then $|W|$ is relatively form 
bounded with respect to $\cE$. Indeed, if $v\in \cD(\cE)$, then
\begin{equation*}\begin{split}
\int_M |W| v^2 d\mu&\le \| W\|_{L^\frac{\nu}{2}} \,\|v^2\|_{L^{\frac{\nu}{\nu-2}}}\\
&\le \| W\|_{L^\frac{\nu}{2}} \left[ A \cE(v)+B \|v\|_{L^2}^2\right]\,. 
\end{split}
\end{equation*}
Moreover, decomposing $|W| = \inf \{ |W|, \lambda \} + W^\lambda$, then for every $v\in \cD(\cE)$,  
\[
\int_M |W| v^2 d\mu\le A \|W^\lambda\|_{L^\frac{\nu}{2}}\, \cE(v)+\left[ \lambda+B\|W^\lambda\|_{L^\frac{\nu}{2}}  \right]  \|v\|_{L^2}^2.
\]
This proves the infinitesimal form boundedness since $\lim_{\lambda\to \infty} \|W^\lambda\|_{L^\frac{\nu}{2}}=0$.
\end{rem}

Another result of the same nature, which is proved exactly as in \cite{ACM}, requires less about $W$ but more 
regularity on the Dirichlet space.
\begin{prop} Let $(M,\mu,\cE)$ be a regular, strongly local Dirichlet space with intrinsic distance compatible with the topology of $M$,
and with Sobolev inequality. Suppose too that the measure $\mu$ is Alfors $\nu$-regular, i.e., there exist constants $0 < c < C$ such that
\[
cr^\nu \le \mu(B(x,r))\le Cr^\nu\quad \mbox{for all}\ x \in M\ \mbox{and}\ r < \mathrm{diam}(M). 
\]
Suppose that $W\in L^q$ for some $q>1$ and moreover, for all $x \in M$ and $r < \mathrm{diam}(M)$, 
\begin{equation}\label{Morrey-1}
\int_{B(x,r)} |W|^q\, d\mu  \leq \Lambda r^{\nu -q\alpha}
\end{equation}
for some constants $\Lambda$ and $\alpha\in [0,2)$. 
If $u\in \cD(\cE)$, $u \geq 0$ and 
\begin{equation*} 
Lu \le Wu\,, 
\end{equation*}
then $u \in L^\infty$.
\end{prop}

It is proved in \cite{ACM} that the Morrey estimate \eqref{Morrey-1} implies that $|W|$ is infinitesimally form bounded with respect to $\cE$. 
In addition, the Gaussian estimate
\[
e^{-tL}(x,y)\le \frac{C}{t^{\nu/2}}\, e^{-\frac{d(x,y)^2}{5t}}\quad {\rm for}\ \ x, y \in M,\ \ t \in (0,1)\,.
\]
is also valid under this hypothesis. 

\subsubsection{Boundedness of solutions of the Yamabe equation}
To apply the results above, we must show that the potential $W$ in \eqref{Yamaberw} satisfies one of these hypotheses.
In fact, any solution to this equation lies in a better $L^p$ space, cf.~\cite{Tr}, \cite{GT}, \cite{Mar} :
\begin{prop}\label{Lp} Let $(M,\mu,\cE)$ be a regular, strongly local Dirichlet space with Sobolev inequality. Suppose that
$W$ is integrable and $W_+$ is infinitesimally form bounded with respect to $\cE$. If $u\in \cD(\cE)$ is a nonnegative solution to
$Lu = Wu$, then $u \in L^q$ for all $q \ge 2$. 
\end{prop}
\proof By assumption on $W_+$, for every $\beta\ge 0$, there are positive constants $A_\beta$ and $B_\beta$ such that 
\begin{equation}\label{Sobolevbeta} 
\|v\|^2_{L^{\frac{2\nu}{\nu-2}}}\le  A_\beta \cE_{-\beta W_+} (v)+B_\beta\|v\|_{L^2}^2\quad \mbox{for every}\ \ v \in \cD(\cE)\,.
\end{equation}

Define, for $\alpha\ge 1$, 
\begin{equation}\label{fa}
f_\alpha(x) = 
\begin{cases}
x^\alpha &{\ \rm if\ } 0\le x\le \alpha^{- \frac{1}{\alpha-1}}\,, \\ 
x + (\alpha^{-\frac{\alpha}{\alpha-1}} - \alpha^{-\frac{1}{\alpha-1}}) &{\ \rm if\ }  \alpha^{-\frac{1}{\alpha-1}}\le x\,. 
\end{cases}
\end{equation}
This function is $\cC^1$ and convex. Next, for $L \geq 1$, set
\[
\phi_{\alpha,L}(x)=L^\alpha f_\alpha\left(\frac{x}{L}\right);
\]
thus $\phi_{\alpha,L}(x) = x^\alpha$ on $[0, \alpha^{-\frac{1}{\alpha - 1}} L]$. If we finally set $G_{\alpha,L}(x)=\int_0^x \phi_{\alpha,L}'(t)^2 \, dt$, 
then a laborious computation gives
\begin{equation}
\phi_{\alpha,L}(x)\le x^\alpha \ \mbox{and}\quad xG_{\alpha,L}(x)\le\frac{\alpha^2}{2\alpha-1}  \left(\phi_{\alpha,L}(x)\right)^2, \quad x \geq 0.
\label{ineqG}
\end{equation}
By the chain rule, with $\varphi=\phi_{\alpha,L}(u)$, 
\begin{multline*}
\cE(\varphi) =\int_M  \phi_{\alpha,L}'(u)^2\, d\gamma(u,u) =\cE(G_{\alpha,L}(u), u) =\int_M G_{\alpha,L}(u)Wu\, d\mu \\
\le \int_M W_+uG_{\alpha,L}(u)\, d\mu \le \frac{\alpha^2}{2\alpha-1} \int_M W_+\varphi^2 \, d\mu\,.
\end{multline*}
Using \eqref{Sobolevbeta} with $\beta=\alpha^2/(2\alpha-1)$ gives
\[
\left\|\phi_{\alpha,L}(u)\right\|^2_{L^{\frac{2\nu}{\nu-2}}}\le B_{\frac{\alpha^2}{2\alpha-1}}\,\, \left\|\phi_{\alpha,L}(u)\right\|_{L^2}^2, 
\]
so that, letting $L\to \infty$, we conclude
\[
u\in L^{2\alpha}\Rightarrow u\in L^{2\frac{\nu}{\nu-2}\alpha}\quad \mbox{for all}\ \ \alpha \geq 1\,.
\]
This completes the proof. 
\endproof

This all leads to the
\begin{prop} Let $(M,\mu,\cE)$ be a regular, strongly local Dirichlet space with Sobolev inequality. Let $V$ be an integrable 
function with nonpositive part $V_-$ infinitesimally form bounded with respect to $\cE$. If $u\in \cD(\cE_V)$ is a 
nonnegative solution to
\[
Lu + Vu = Y(V) u^{\frac{\nu+2}{\nu-2}}\,,
\]
then for every $p \geq 2$, 
\[
\int_M u^p \, d\mu<\infty .
\]
\end{prop}
Indeed, the assumption that $u\in \cD(\cE_V)$ and \eqref{eq:Sobolev} give that $u^{\frac{4}{\nu-2}} \in L^{\frac{\nu}{2}}$. 
According to Remark~\ref{nu2}, $u^{\frac{4}{\nu-2}}$ is infinitesimally form bounded with respect to $\cE$. 
We can thus apply Proposition~\ref{Lp} with $W = - V + Y(V) u^{\frac{4}{\nu-2}}$. 

\subsection{Positivity of solutions}
The argument of \cite{ACM} , see also \cite{Gursky}, can be applied verbatim to our Yamabe equation 
when $(M,\mu,\cE)$ is a regular, strongly local Dirichlet space with intrinsic distance compatible with the topology of $M$. 
Thus any nonnegative solution of this equation which is strictly positive on some ball is strictly positive everywhere, 
provided that \eqref{eq:Sobolev} holds and $|V|$ satisfies a Morrey type estimate. However, the Harnack
inequality need not hold in this generality. In the next subsection, we give a criterion which ensures H\"older continuity 
of solutions to the linear equation $Lu = f$, and this implies that if $u \not\equiv 0$, then it is 
strictly positive on some ball. 

\subsection{Higher regularity of solutions}
We now turn to questions concerning the modulus of continuity of solutions of the equation $Lu = f$. 
As usual, let $(M,\mu,\cE)$ be a regular, strongly local Dirichlet space with intrinsic distance compatible with the topology of $M$. 
We assume that the measure $\mu$ is Ahlfors $\nu$-regular.
and a uniform Poincar\'e inequality holds. This means that if $r \leq \frac{1}{4} \diam M$, then
\begin{equation*}
\left\|v-v_B\right\|^2_{L^2(B)}\le C r^2 \int_{B(x,2r)} d\gamma(v,v)\quad \mbox{for every}\ \  v \in \cD(\cE)\,, 
\end{equation*} 
where $B = B(x, r)$ and $v_B=\frac{1}{\mu(B)} \int_B v \, d\mu$. 
For a nice review on the Dirichlet space satisfying these assumptions, see \cite{SC2} and also the paper \cite{GriTelcs} for recent results.

These assumptions imply that the heat kernel of $L$ exists and satisfies 
Gaussian upper bounds, and also that the Sobolev inequality \eqref{eq:Sobolev} holds.  They also guarantee the
elliptic and parabolic Harnack inequality. In particular, if $h$ is a positive harmonic function on $2B := B(x, 2r)$ (so $Lh = 0$ on 
$2B$), then
\[
\sup_{z\in B} h(z) \le C_\cH \inf_{z\in B} h(z)\,.
\]
The Harnack constant $C_\cH$ depends only on the constants in the Ahlfors regularity condition and the Poincar\'e inequality.
From this, one obtains H\"older continuity of harmonic functions.
\begin{lem} 
Let $h \in L^\infty(2B)$ be a solution of the equation $Lh=0$ on a ball $2B$. Then, for all $p, q \in B(x,r)$, 
\[
\left| h(p)-h(q)\right|\le C \left(\frac{d(p,q)}{r}\right)^\beta \sup_{z\in 2B}|h(z)|\,.
\]
\end{lem}
In fact, $\beta=\log_2\left(\frac{C_\cH+1}{C_\cH-1}\right)$ and $C=2 (\frac{C_\cH+1}{C_\cH-1})$.
The Green function of $L$ is a symmetric function $G\in \cC^0((M\times M) \setminus \mathrm{Diag})$ 
such that $G(x, \cdot )\in L^1(M)$ for any $x \in M$, and in addition, if $f \in L^2$, then 
\[
\cD(\cE) \ni u(x) = \int_M G(x,y)f(y)\, d\mu(y)\quad  \mbox{for $\mu$-a.e.} \  x \in M 
\]
satisfies $Lu = f - f_M$, where $f_M=\frac{1}{\mu(M)}\int_M f\, d\mu$. Clearly, if $Lu=f$ and $\int_M u \, d\mu=0$, 
then $f_M = 0$ and 
\[
u(x) = \int_M G(x,y)f(y)\, d\mu(y)\quad \mbox{for $\mu$-a.e.} \  x \in M. 
\]

\begin{prop}
The Green kernel $G$ satisfies 
\[ 
|G(x,y)|\le \frac{C}{d(x,y)^{\nu-2}}\quad \mbox{for all}\  x, y \in M 
\]
and if $p, q, y \in M$ with $d(p,q) \le \frac{1}{2} d(p,y)$, then 
\[
\left|G(p,y)-G(q,y)\right | \le C  \left(\frac{d(p,q)}{d(p,y)}\right)^\beta  \frac{1}{d(q,y)^{\nu-2}}\,.
\]
\end{prop}
\begin{thm} 
With all the assumptions as above, suppose that $f$ satisfies the Morrey estimate
\begin{equation*}
\int_{B(x,r)} |f|\, d\mu  \leq \Lambda r^{\nu -\alpha}\quad \mbox{for all}\  x \in M,\ r \leq \frac{1}{2}\diam M\,,
\end{equation*}
for some $\Lambda > 0$ and $\alpha\in [0,2)$. If $u\in \cD(\cE)$ solves $Lu = f$, then $u$ is H\"older 
continuous of order $\mu = \min\{\beta, 2-\alpha\}$ when $\beta \not= 2 \alpha$; if $\beta = 2\alpha$,
then $\mu$ need only satisfy $0 < \mu < \min\{\beta, 2-\alpha\}$.
\end{thm}
\proof
For each $x\in M$, introduce the nondecreasing function 
\[
\nu_x(r)=\int_{B(x,r)} |f|d\mu\,.
\]
If $p, q \in M$ and $\rho:=d(p,q) \le \frac{1}{8}\diam M$, then 
\begin{equation*}\begin{split}
\left|u(p)-u(q)\right|\le &\int_{B(p,4\rho)} |G(p,y)| |f(y)|\, d\mu(y)+\int_{B(q,4\rho)} |G(q,y)| |f(y)| \, d\mu(y)\\
&+\int_{M\setminus B(p,4\rho)} \left|G(p,y)-G(q,y)\right|\,  |f(y)|\, d\mu(y)\,. 
\end{split}\end{equation*}
Using the estimates on $G$ and integrating by parts, 
\begin{multline*}
\int_{B(p,4\rho)} |G(p,y)| |f(y)|\, d\mu(y)\le \\
\int_0^{4\rho} \frac{C}{r^{\nu-2}} \, d\nu_p(r) 
=C\left( \frac{\nu_p(4\rho)}{4^{\nu-2}\rho^{\nu-2}}+(\nu-2)\int_0^{4\rho}\frac{\nu_p(r)}{r^{\nu-1}}\, dr \right) 
\end{multline*}
so by the Morrey estimate, 
\[
\int_{B(p,4\rho)} |G(p,y)| |f(y)|\, d\mu(y)\le C \rho^{2-\alpha}.
\]
The integral over $B(q,4\rho)$ is bounded by $C \rho^{2-\alpha}$ too. 

For the final term,
\begin{multline*}
\int_{M\setminus B(p, 4\rho)} \left|G(p,y)-G(q,y)\right|\,  |f(y)|\, d\mu(y) \\ 
\le C\rho^\beta \, \left[ \frac{\nu_p(\diam M)}{(\diam M)^{\nu-2+\beta}} +
(\nu-2+\beta)\int_{4\rho}^{\diam M} \frac{\nu_p(r)}{r^{\nu-2+\beta}} \, dr \right].
\end{multline*}

Since 
\[
\int_{4\rho}^{\diam M} \frac{\nu_p(r)}{r^{\nu-2+\beta}} \, dr\le\Lambda  \int_{4\rho}^{\diam M} \frac{1}{r^{\alpha-2+\beta}} \, dr
\]
and
\[
 \int_{4\rho}^{\diam M} \frac{1}{r^{\alpha-1+\beta}} dr\le \begin{cases}
 \frac{1}{\beta-2+\alpha} \frac{1}{(4\rho)^{\alpha-2+\beta}} &\,\mathrm{if}\ \ 2 - \alpha < \beta\,, \\[3mm]
  \frac{1}{\alpha+2-\beta} \frac{1}{(\diam M)^{\alpha-2+\beta}}& \,\mathrm{if}\ \ \beta < 2 - \alpha\,, \\[3mm]
 \log( \diam M/ (4\rho)) & \,\mathrm{if}\ \ \beta = 2 - \alpha\,, 
\end{cases} 
\]
we conclude that $u$ is H\"older continuous of order $\mu$. 
\endproof

\subsection{Conclusion}
\begin{thm}
Let $(M,\mu,\cE)$ be a regular, strongly local Dirichlet space with Sobolev inequality.
\begin{enumerate}[$i)$]
\item If $V_-\in L^q$ for some $q>\nu/2$ and if $Y(V) < 1/A$, then there exists a nonnegative bounded function 
$u \in \cD(\cE)$ such that 
\[
\cE_V(u)=Y(V)\,\, \mbox{and}\,\, \|u\|_{L^{\frac{2\nu}{\nu-2}}}=1\,.
\]
\item The conclusion of i) holds if the intrinsic distance is compatible with the topology of $M$, if $\mu$ is Ahlfors $\nu$-regular, 
and if $V_-$ satisfies the Morrey estimate \eqref{Morrey-1}.
\item Suppose that the intrinsic distance is compatible with the topology of $M$, $\mu$ is Ahlfors $\nu$-regular,
and $|V|$ satisfies the Morrey inequality \eqref{Morrey-1}. (Note that this holds if $V \in L^p$ for some $p > \nu/2$.) 
Assume also that $(M, \mu, \cE)$ satisfies the Poincar\'e inequality on any ball $B = B(x,r) \subset M$, $r \leq \frac{1}{4} \diam M$,
and that $Y(V) < 1/A$. 
Then the solution in either case i) or case ii) is strictly positive and H\"older continuous.
\end{enumerate}
\end{thm}

\section{Examples}
We now explain how the general results above simplify the original proof of the usual Yamabe problem and then 
yield a generalization of the CR Yamabe problem. 

\subsection{The Riemannian Yamabe problem} 
We have already discussed the classical Yamabe problem when $M^n$ is a compact smooth Riemannian manifold,
$n \geq 3$. The metric $\widehat g=v^{\frac{4}{n-2}} g$ has constant scalar curvature if and only if $v$ is a critical point 
of the Yamabe functional
\begin{equation*}\begin{split}
Q_g(f)&:=\frac{\int_M \left[\frac{4(n-1)}{n-2}|df|^2_g + \scal_g f^2\right] d\mu_g}{\left(\int_M  f^{\frac{2n}{n-2}} d\mu_g \right)^{1-\frac{2}{n}}}\\
&={\rm Vol}_{\widetilde{g}}(M)^{\frac{2}{n}-1}\int_M\scal_{\widetilde{g}}d\mu_{\widetilde{g}}\quad \mathrm{for}\ \ \widetilde{g}=f^{\frac{4}{n-2}} g. 
\end{split}\end{equation*}
The minimizer for this problem always exists. For this case, 
\begin{itemize}
\item the pair of $(M, d\mu_g)$ and $\mathcal{E}(v) = \frac{4(n-1)}{n-2}\int_M |dv|^2_g \, d\mu_g$, $v \in W^{1,2}(M, d\mu_g)$,
determine the Dirichlet space; 
\item $L + V = -\frac{4(n-1)}{n-2}\Delta_g + {\rm Scal}_g$; 
\item $\nu = n$ and $Y_{\ell} = Y(S^n, [g_0])$.
\end{itemize}
(Here $\Delta_g = \mathrm{div} \nabla$.) 

The key result, due to Aubin~\cite{Aubin} and Schoen~\cite{Sc-1}, \cite{Sc-2}, \cite{SY}, states that if $(M, [g])$ is 
not conformal to $(S^n, [g_0])$, then $Y(M, [g]) < Y(S^n,[g_0])$, so the existence proof above may be applied. 

\subsection{The contact Riemannian Yamabe problem} 
The second application of our results is to the Yamabe problem on contact Riemannian manifold.  This problem
was initially posed for CR manifolds by Jerison and Lee, \cite{JerisonLee1}, \cite{JerisonLee2}, and solved by them 
for manifolds not CR equivalent to the standard sphere $S^{2m +1}$. The remaining case was completed by Z. Li \cite{Liz}. 
This problem does not seem to have been treated for non-integrable almost complex structures, but we
are able to work in that more general context here.

\subsubsection{The setting:} Recall from \cite{DT}, \cite{Tanno} that a {\it contact manifold} is an odd dimensional manifold 
$M^{2m+1}$ with a totally non-integrable hyperplane subbundle $H \subset TM$. Thus, for each $x \in M$, there
is a nondegenerate bilinear form
\begin{equation*}
H_x \times H_x\rightarrow T_xM/H_x,\ \ (X,Y)\mapsto [X,Y] \ \ \mbox{mod}\ H_x.
\end{equation*}
When $M$ and $H$ are oriented, one may choose a {\it contact form} $\theta \in \Omega^1(M)$ with $H=\ker \theta$;
in terms of this form, nondegeneracy of $H$ is equivalent to $\theta\wedge (d\theta)^m\neq 0$ everywhere on $M$. 

A choice of $\theta$ uniquely defines the {\it Reeb vector field} $\xi \in \frak X (M)$; this is associated to $\theta$ by the conditions
\[
\theta(\xi) = 1\ \ {\rm and}\ \ \cL_\xi \theta = 0,
\]
(here $\cL_{\xi}$ is the Lie derivative of $\xi$). Thus 
\begin{equation}\label{decomp} 
TM = H \oplus \mathbb{R}~\xi\,. 
\end{equation} 

A {\it contact Riemannian manifold} $(\theta, g_H, J)$ is a triple, consisting of a contact form $\theta$, a Riemannian metric $g_H$ 
on $H$, and a compatible almost complex structure $J$ on $H$, i.e.\ such that 
$$ 
g_H(X, Y) = d\theta(JX, Y)\quad \mbox{for every}\ \ x \in M,\ \ X, Y \in H_x\,. 
$$ 
For any contact form $\theta$, there always exists a compatible pair $(g_H, J)$, see \cite{Blair}.  We can then define
the {\it Webster metric} $g_{\theta}$ on $M$, which is Riemannian, by
\[
g_{\theta} = \pi_H^{\ast}g_H + \theta^2, 
\]
where $\pi_H : TM \rightarrow H$ is the projection associated with the decomposition \eqref{decomp}. By definition,
$\xi \perp H$ with respect to $g_{\theta}$. We also define the {\it Tanaka-Webster scalar curvature} $\scal_{g_H}$ by
\[
\scal_{g_H} = \scal_{g_\theta} - {\rm Ric}_{g_\theta}\left(\xi,\xi\right) + 4m\, . 
\]

The structure $(\sigma \theta, \sigma g_H, J)$, where $\sigma \in \calC^\infty(M)$, $\sigma > 0$, is said to be
conformally related to $(\theta, g_H, J)$, and the conformal class $[\theta, g_H, J]$ is the set of all such conformally
related structures. Associated to $(\sigma \theta, \sigma g_H, J)$ are its Reeb vector field
\[
\xi_{\sigma \theta} = \frac{1}{\sigma} \Big{(} \xi_\theta + \frac{1}{2 \sigma} J \pi_H(\nabla^{g_\theta} \sigma) \Big{)} 
\]
and Webster metric 
\[ 
g_{\sigma \theta} = \sigma (\pi_H^{\ast} g_H) + \sigma^2 \theta^2\, .
\]
The contact Riemannian analogue of the conformal Laplacian is the operator 
\[
- b_m \Delta_H u + \scal_{g_H} u, \qquad b_m = \frac{4(m+1)}{m},
\]
where
\[
\Delta_H = \Delta_{g_{\theta}} - \xi \circ \xi 
\]
is the horizontal Laplacian, which is defined as follows.  For any function $v$, consider the restriction of $dv$ to $H$.
This has squared $g_H$-length
\[
|d_Hv|^2(x) = \sup_{X\in H_x \atop g_H(X,X)\le 1} |dv(X)|^2 = |dv|_{g_{\theta}}^2(x) - |\xi v|^2(x)\, .
\]
Integrating this against the volume form $\theta\wedge(d\theta)^m$ gives a quadratic form, and $\Delta_H$ is then
determined by 
\[
\int_M |d_Hv|^2\,  \theta\wedge(d\theta)^m = - \int_M v\Delta_Hv\, \theta\wedge(d\theta)^m,\  v \in \calC^\infty(M)\, . 
\]

The distance function $\rho$ associated to this quadratic form is the sub-Riemannian distance
\begin{multline*}
\rho(x,y) = \inf\big\{\int_0^1 |\dot c|^2_{g_H}\, dt:  c \in \calC^1([0,1], M), \\
c(0)=x, c(1)=y, \ \ \dot c(t) \in H_{c(t)}\ \mbox{for all}\  t \big\}, 
\end{multline*}
If $d\gamma$ is the energy measure associated to the quadratic form, then 
\[
\rho(x,y) = \sup\left \{u(x)-u(y) : u \in \Lip(M),\ \ d\gamma(u, u) \leq \theta \wedge (d\theta)^m \right\}.
\]
Note that $\rho$ is compatible with the geodesic distance for $d_{g_{\theta}}$ in the sense that 
\[
d_{g_{\theta}} \leq \rho\le C \sqrt{d_{g_{\theta}}}
\]
for some $C > 0$. 

Just as in Riemannian geometry, there is a simple conformal transformation rule for the Tanaka-Webster scalar curvature. 
If $\sigma = u^{\frac{4}{\alpha - 2}}$, $\alpha = 2m + 2$, then writing $\widehat \theta = \sigma \theta$, 
$\widehat g_H = \sigma g_H$, we have
\begin{equation}\label{CTR}
- b_m \Delta_H u + \scal_{g_H} u = \scal_{\widehat g_H} u^{\frac{\alpha + 2}{\alpha - 2}} \, .
\end{equation}
Noting that 
\[
\int_M \scal_{\widehat g_H}\, \widehat{\theta}\wedge(d\widehat{\theta})^m = 
\int_M\left[b_m |d_Hu|^2+\scal_{g_H} u^2\right]\, \theta\wedge(d\theta)^m
\]
and 
\[
\int_M \, \widehat{\theta}\wedge(d\widehat{\theta})^m = \int_M |u|^{\frac{2\alpha}{\alpha-2}}\, \theta\wedge(d\theta)^m,
\]
we define the {\it contact Yamabe invariant} of $(M, [\theta, g_H])$ by
\[
Y\left(M,\left[\theta,g_H\right]\right) = 
\inf_{u  > 0} \frac{\int_M\left[b_m |d_Hu|^2+\scal_{g_H} u^2\right]\, \theta\wedge(d\theta)^m}
{\left(\int_M |u|^{\frac{2\alpha}{\alpha-2}}\, \theta\wedge(d\theta)^m\right)^{1-\frac{2}{\alpha}}}.
\]
If this infimum is attained by some $(\widehat \theta,\widehat g_H)$, then the Euler-Lagrange equation for this 
functional shows that $(M, \widehat \theta,\widehat g_H)$ has constant Tanaka-Webster scalar curvature. 

\subsubsection{The Heisenberg group:} 
The basic model contact Riemannian manifold is the {\it Heisenberg group} 
\[
\mathfrak{h}_m = (\R^{2m+1},\theta_0):=\left(\left\{(x,y,t)\in \R^m\times\R^m\times\R \right\},\ \theta_0 = dt - \sum_j y_jdx_j\right) 
\]
with metric $g_0 = |dx|^2 + |dy|^2$ on the horizontal distribution 
\[
H = \ker \theta_0 = {\rm span}_{\mathbb{R}} \Big{\{} \frac{\partial}{\partial x_j}-y_j\frac{\partial}{\partial t}, \frac{\partial}{\partial y_j},\ \ 
j = 1, \cdots, m \Big{\}}. 
\]
It is simple to check that $\scal_{g_0} = 0$, while 
\[
Y(\mathfrak{h}_m, [\theta_0, g_0]) > 0. 
\]
(Since $\mathfrak{h}_m$ is noncompact, this invariant is the infimum over compactly supported smooth nonnegative functions.)

The Heisenberg group satisfies a uniform Poincar\'e-Wirtinger inequality: there exists a constant $C > 0$ such that 
on any $\rho$-ball $B_r$ of radius $r > 0$,
\[
\int_{B_r} u^2~\theta_0 \wedge(d\theta_0)^m \leq C r^2 \int_{B_r}  |d_Hu|^2\, \theta_0 \wedge (d\theta_0)^m 
\]
for all $u \in L^1(B_r)$ with
\[
|d_Hu|^2 \in L^1(B_r)\quad {\rm and}\quad \int_{B_r} u~\theta_0 \wedge (d\theta_0)^m = 0. 
\]

\subsubsection{The local Yamabe invariant:} 
Jerison and Lee \cite{JerisonLee1} posed the problem of showing that on a given CR manifold, 
the local Yamabe invariant equals the contact Yamabe invariant of the Heisenberg group in normal CR-coordinates. 
It turns out that if one uses Darboux coordinates instead, then this is not difficult. Indeed, let $(M^{2m+1}, \theta, g_H)$ 
be a compact contact Riemannian manifold. For each $p \in M$, there exists a diffeomorphism 
$$ 
\varphi \colon U\rightarrow \bB(1) = \{ (x, y, t) \in \R^{2m+1} : \rho( (x, y, t), {\bf 0}) < 1 \},
$$ 
where $U$ is a neighbourhood of $p$, such that 
\[
\varphi^*\theta_0=\theta,\quad (\varphi^*g_0)_p = (g_H)_p\quad \mathrm{and}\quad 
\varphi^*g_0-g_H = O(\sqrt{\varepsilon})\ \ {\rm on}\ \ \cU_\eps := \varphi^{-1}(\bB(\eps))\,.  
\]
Note that the Tanaka-Webtser scalar curvature is assumed to be bounded, hence the local Yamabe invariant is equal to the local Sobolev invariant.
The Heisenberg group has vanishing Tanaka-Webster scalar curvature 
and the Sobolev constant varies continuously when the metric varies in the space of continuous metrics, 
hence, we have 
\begin{equation*}\begin{split}
Y\left(\cU_\eps,\left[\theta,g_H\right]\right)&=Y\left(\bB(\eps),\left[\theta_0,g_0\right]\right)\,\left(1+O\left(\sqrt{\varepsilon}\right)\right)\\
&= Y(\mathfrak{h}_m, [\theta_0, g_0])\,\left(1+O\left(\sqrt{\varepsilon}\right)\right)\,.
\end{split}\end{equation*}

Since the Heisenberg group satisfies the Poincar\'e-Wirtinger inequality,  
it is easy to show that any compact contact Riemannian manifold $(M^{2n+1},\theta,g_H)$ 
satisfies a local Poincar\'e-Wirtinger inequality.  There exist positive constants $C, R > 0$ such that 
on any ball $B_r \subset M$ of radius $r \in (0,R)$, 
$$ 
\int_{B_r} u^2~\theta \wedge (d\theta)^m \leq C r^2 \int_{B_{2r}}  |d_Hu|^2\, \theta \wedge (d\theta)^m 
$$ 
for any $u \in L^1(M)$ with 
\[
|d_Hu|^2 \in L^1(B_{2r})\quad {\rm and}\quad \int_{B_r} u~\theta \wedge (d\theta)^m = 0\,. 
\]

\subsubsection{The contact Riemannian Yamabe problem:} 
The second application of our general result is to the contact Riemannian Yamabe problem. This is
a generalization of the first main result of Jerison and Lee.
\begin{thm}Let $(M^{2m+1}, \theta,g_H)$ be a compact contact Riemannian manifold. 
Assume that 
\[
Y(M, [\theta, g_H]) < Y(\mathfrak{h}_m, [\theta_0, g_0])\,.  
\]
Then there exists a positive function $u\in \calC^\infty(M)$ such that the Tanaka-Webster scalar curvature of
$(u^{\frac{2}{m}}~\theta,\ u^{\frac{2}{m}}~g_H)$ is constant. 
\end{thm} 

In this setting, 
\begin{enumerate}
\item[$\bullet$] $d\mu = \theta \wedge (d\theta)^m$, 
\item[$\bullet$] $\mathcal{E}(v) = \frac{4(m+1)}{m}\int_M |d_Hv|^2~\theta \wedge (d\theta)^m$ 
\item[$\bullet$] $v \in \cD(\cE) = \{ v \in L^2(M): |d_Hv|^2 \in L^1(M) \}$, 
\item[$\bullet$] $L + V = - \frac{4(m+1)}{m} \Delta_H + {\rm Scal}_{g_H}$, 
\item[$\bullet$] $\nu = 2m + 2$,
\item[$\bullet$] $Y_{\ell} = Y(\mathfrak{h}_m, [\theta_0, g_0])$. 
\end{enumerate}

Although our result gives only a positive bounded solution $u \in \cD(\cE)$, the hypoelliptic properties of 
$\Delta_H$ directly show that $u \in \calC^\infty$.

\end{document}